\documentclass[10pt]{amsart}
\def\dim{\operatorname{dim}}
\def\Ker{\operatorname{Ker}}
\newcommand{\ZZ}{\mathbb Z}
\newcommand{\NN}{\mathbb N}
 \def\nt{\noindent} 
\def\ini{\mathop{\kern0pt\fam0in}\nolimits}
\def\Rees{{\mathcal R}}
\def\rk{\operatorname{rk}}
\def\cocoa {\mbox{\rm C\kern-.13em o\kern-.07 em C\kern-.13em o\kern-.15em A}}
\let\to=\rightarrow

\newtheorem{theorem}{Theorem}[section]
\newtheorem{lemma}[theorem]{Lemma}
\newtheorem{corollary}[theorem]{Corollary}
\newtheorem{proposition}[theorem]{Proposition}
\newtheorem{remark}[theorem]{Remark}
\newtheorem{example}[theorem]{Example}
\newtheorem{definition}[theorem]{Definition}
\newtheorem{conjecture}[theorem]{Conjecture}

\def\addots{\mathinner{\mkern1mu\raise1pt\hbox{.}\mkern2mu\raise4pt\hbox{.}
         \mkern2mu\raise7pt\vbox{\kern7pt\hbox{.}}\mkern1mu}}

\begin{document}
\title{Linear spaces, transversal polymatroids and ASL domains}
\author{Aldo Conca  }
\address{Dipartimento di Matematica, Universit\'a di Genova,  Genova, Italy}  
\email{conca@dima.unige.it}
  \maketitle

\section{Introduction}
\label{intro}
 Let $K$ be an infinite field and $R=K[x_1,\dots,x_n]$ be the   polynomial ring.  Let $V=V_1, \dots, V_m$ be a
collection of vector spaces of linear forms. Denote by
$A(V)$ the $K$-subalgebra of $R$ generated by the elements of the product $V_1\cdots V_m$. Our goal is to investigate  the
properties of the algebra
$A(V)$ and the  relations with two  problems in algebraic combinatorics:  White's and related  conjectures on polymatroids and the
study of integral posets. 

\subsection*{ Polymatroids}  A finite subset  $B$ of $\NN^n$  is the base set  of a discrete  polymatroid $P$  if for every
$v=(v_1,\dots,v_n), w=(w_1,\dots,w_n)\in B$ one has $v_1+\dots+v_n=w_1+\dots+w_n$ and for all $i$ such that $v_i>w_i$ there
exists a $j$ with
$v_j<w_j$ and $v+e_j-e_i\in B$. Here $e_k$ denotes the $k$-th  vector of
the standard basis of $\NN^n$. The notion of discrete polymatroid is a generalization of the
classical notion of matroid, see \cite{Ed, HH, O,Wel}. Associated with the base
$B$ of a discrete polymatroid
$P$ one has a $K$-algebra $K[B]$,   called  the base ring of $P$,  defined to be the $K$-subalgebra of $R$ generated by the monomials
$x^v$ with $v\in B$. The algebra $K[B]$ is known to be normal and  hence Cohen-Macaulay  \cite{HH}.  White predicted in
\cite{Wh} the shape of the defining equations of $K[B]$ as a quotient of a polynomial ring: they should the quadrics arising from  the
so-called symmetric exchange relations of the polymatroids.  Herzog and Hibi 
\cite{HH}  did not ``escape from the temptation'' to ask  whether  $K[B]$ is  defined by a Gr\"obner basis of
quadrics  and whether  $K[B]$  is  a Koszul algebra. These two questions are closely related to  White's conjecture. This is  because
for any standard graded algebra $A$ with defining ideal $I$, the existence of a Gr\"obner basis of quadrics for $I$   implies
the Koszul property  of $A$ which implies that $I$ is defined by quadrics. 

If   $C_1,\dots,C_m$   are  non-empty subsets  of
$\{1,\dots,n\}$  then  the set of the vectors  $\sum_{k=1}^m  e_{j_k}$  with $j_k\in C_k$ is the base of a polymatroid.  Polymatroids
of this kind are called transversal.  Therefore the  base rings of transversal polymatroids are exactly the rings of type $A(V)$ where
the spaces $V_i$ are generated by variables. 
 For  transversal polymatroids we prove that the base ring  $K[B]$ is  Koszul   and describe the  defining equations, see Section
 \ref{mono}. Indeed,   $K[B]$ is defined as a quotient of a Segre product $T^*$ of polynomial rings   by a Gr\"obner basis of linear
binomial forms  of $T^*$. 

\subsection*{ASL and integral posets}  Algebras with straightening laws   (ASL for short) on  posets  were introduced by  De Concini,
Eisenbud and Procesi   \cite{DEP, Ei}, see also \cite{BV}.   The  abstract definition of ASL was   inspired by earlier   work of 
  Hochster,  Hodge,   Laksov, Musili, Rota,  and Seshadri among others.  It was  motivated by the existence of many 
families of  classical algebras,   such as  coordinate rings of Grassmannians and their  Schubert subvarieties  and
various kinds of determinantal  rings, which could be treated within that  frame.     We recall in   \ref{ASL} the definition  of
homogeneous ASL and in \ref{ASL-GB} a well-known characterization of them  in terms of  revlex  Gr\"obner bases.  

A finite  poset $H$ is   integral (with respect to a field $K$) if there exists a homogeneous ASL  
domain supported on $H$.    A beautiful  result,  due to Hibi 
\cite{H1},  says  that any distributive lattice $L$ is integral. Indeed, $L$  supports a homogeneous  ASL domain, denoted by 
$H_L$, in a very natural way.  The ring
$H_L$ is called the Hibi ring of $L$ and its  defining equations  are the so-called Hibi relations:   
$xy-(x\wedge y)(x\vee y)$.  In a series of papers \cite{HW1,HW2,H2,W1,W2}   Hibi and Watanabe classified   various families  of
integral posets of low dimension.    In this direction,  we construct a new class of integral posets: the rank truncations of
hypercubes. In details, given a sequence of positive integers
$d=d_1,\dots,d_m,$   let
$H(d)=\Pi_{i=1}^m\{1,\dots, d_i\}$ and, for $n\in \NN$,  $H_n(d)=\{
\alpha\in H(d) : \rk \alpha< n\}$.  We show  that $H_n(d)$ is an integral poset (over every infinite field
$K$). This is done by proving that $A(V)$ is a homogeneous  ASL on $H_n(d)$ if the $V_i$ are generic linear spaces of dimension $d_i$
of $R$, see  Section \ref{gene}.  In particular, our construction shows  that the Veronese subrings of polynomials rings are
homogeneous ASL (obviously domains). Note however that they are not,  in general,  ASL  with respect to  their  semigroup
presentation. 

\subsection*{}  Results from \cite{CH} show that  for any collection $V=V_1,\dots,V_m$  
the algebra $A(V)$ is normal.  As said above, in the monomial case,  i.e. when the $V_i$ are generated by variables, 
  we show that $A(V)$ is Koszul and describe its defining equations.      Our argument for the monomial case is based on a certain
elimination   process and on a result,  Theorem \ref{StuVil},  proved independently by Sturmfels and Villarreal,  describing the
universal Gr\"obner basis of the ideal of
$2$-minors of a matrix of variables.    This approach  suggests also  a possible strategy for proving that $A(V)$ is Koszul in  the
general case. The elimination   process is still available and what one needs is a replacement of the Sturmfels-Villarreal's theorem. 
This boils down to the following:

\begin{conjecture}
\label{strocon} 
Let $t_{ij}$  be  distinct variables over a field $K$ with $1\leq i\leq m$ and $1\leq j\leq n$.  Let
$L=(L_{ij})$ be a
$m\times n$ matrix with $L_{ij}=\sum_{k=1}^n a_{ijk}t_{ik}$ and
$a_{ijk}\in K$ for all $i,j,k$.  Denote by $I_2(L)$ the ideal of the
$2$-minors of $L$.  We conjecture that for  every choice  of $a_{ijk}$'s,  for every  term order $<$ on $K[t_{ij}]$ the   initial ideal
$\ini_<(I_2(L))$  is  square-free in the $\ZZ^m$-graded sense, i.e. it is generated by elements  the form  $t_{i_1j_1}\cdots
t_{i_kj_k}$ with $i_1<i_2<\cdots <i_k$. 
\end{conjecture} 

This conjecture can be rephrased    in terms of universal comprehensive Gr\"obner bases \cite{We}: the parametric ideal $I_2(L)$
(the parameters being the $a_{ijk}$'s)  has a comprehensive and universal Gr\"obner basis whose elements are multihomogeneous of
degree bounded by
$(1,1,\dots,1)$. 
 
If  $L=(t_{ij})$    then   \ref{strocon} holds; this is a consequence of Theorem \ref{StuVil}. We prove in
\ref{main} that \ref{strocon} holds  when the $a_{ijk}$ are  generic.  As a consequence, we are able to show that for generic spaces
$V_i$    the algebra $A(V)$ is Cohen-Macaulay, Koszul and describe the defining equations of $A(V)$. 
 In particular, as mentioned above,  in the generic case $A(V)$ turns out to be a homogeneous  ASL on the poset  $H_n(d)$ where
$d=d_1,\dots,d_m$ and $d_i=\dim V_i$.  

We thank   C.Krattenthaler  who provided  a combinatorial argument for a statement which was used in an earlier version of   
the proof of \ref{main}.  The  results   presented in  this  the paper  have been inspired, suggested and confirmed by  computations  
performed  by   the computer algebra system CoCoA  
\cite{Cocoa}.

\section{Normality of $A(V)$}
\label{norma}

Let $I_i$ be the ideal of $R$ generated by $V_i$.  In  \cite{CH}  it is proved that the product ideal
$I_1\cdots I_m$ has always a linear resolution. One of the main step in proving that result is the following \cite[3.2]{CH}:

\begin{proposition} 
\label{primdec} 
For any subset $A\subseteq  \{1,\dots, m\}$ set $I_A=\sum_{i\in A} I_i$ and denote by
$\#A$ the cardinality of $A$. Then 
$$I_1\cdots I_m = \cap  I_A^{\#A}$$ is a primary decomposition of $I$. Here the intersection is extended to all the $A\neq
\emptyset$. 
\end{proposition} 

Proposition \ref{primdec} easily implies:

\begin{theorem} 
\label{normal}
$A(V)$ is normal.
\end{theorem} 

\begin{proof}  Set $J=I_1\cdots I_m$.  Note that   $I_A$ is a prime ideal generated by linear forms. Hence the powers of $I_A$ are
integrally closed. It follows that $J$ is integrally closed. Since the powers of $J$ are again product of ideals of linear forms, the
same argument  apply also to the powers of $J$.  Hence we   conclude that $J$ is normal (i.e. all the powers of $J$ are integrally
closed). This is equivalent to the fact that the Rees algebra $\Rees(J)=\oplus_{k\in
\NN}  J^k$ is normal. Now $A(V)$, being a direct summand of $\Rees(J)$, is normal as well.
\end{proof}

\section{The monomial case}
\label{mono}

We  now analyze the monomial case.  Our goal is to show that $A(V)$ is Koszul if each $V_i$ is monomial and to develop a strategy to
attack the general case. So in this section we assume that each
$V_i$ is generated by a subset of the variables $\{x_1,\dots,x_n\}$.  Say $V_i=\langle x_j : j\in C_i
\rangle$ where $C_i$ is a non-empty  subset of $\{1,\dots, n\}$.  Consider the auxiliary algebra 

$$B(V)=K[V_1y_1,\dots, V_my_m]=K[y_ix_j : i\in 1,\dots, n, \mbox { and } j\in C_i]$$

where $y_1,\dots, y_m$ are new variables.  The algebra $B(V)$ sits inside the Segre product 

$$S=K[y_ix_j : 1\leq i\leq m, 1\leq j\leq n].$$  

We consider variables $t_{ij}$ with $i=1,\dots,m$ and $j=1,\dots,n$, and define 

$$T=K[t_{ij}  : 1\leq i\leq m, 1\leq j\leq n]    \quad \quad \mbox{ and }  \quad \quad  T(V)=K[t_{ij}  : 1\leq i\leq m,j\in C_i]$$

and presentations: 

$$\phi:T\to S    \quad   \mbox{ and }  \quad  \phi':T(V)\to B(V)  $$ defined by sending $t_{ij}$ to $y_ix_j$.

It is well-known  that $\Ker \phi$ is the ideal $I_2(t)$ of  $2$-minors of the $m\times n$ matrix $t=(t_{ij})$.   Then the algebra
$B(V)$ is defined as a quotient of $T(V)$ by the ideal  $I_2(t)\cap T(V)$.     The algebras  $B(V), T(V), S$ and $T$  can be given a 
$\ZZ^m$-graded structure by setting the degree of
$y_ix_j$ and
$t_{ij}$ to be $e_i\in \ZZ^m$.

By work of Sturmfels \cite[4.11 and  8.11]{S} and Villarreal \cite[ 8.1.10 ]{V}  one knows  that
 a universal Gr\"obner basis of $I_2(t)$ is given by the  cycles of the complete bipartite  graph $K_{n,m}$. In details, a 
cycle of the complete bipartite graph is  described  by a  pair 
$(I,J)$ of sequences of   integers, say 

$$I=i_1,\dots,i_s, \quad J=j_1,\dots, j_s$$ with  $2\leq s\leq \min(n,m)$,  $1\leq i_k\leq m$,  $1\leq j_k\leq n$,  and such that the
$i_k$ are distinct and the
$j_k$ are distinct.   Associated with   any such a pair we  have   a polynomial 
$$f_{I,J}=t_{i_1j_1}\cdots t_{i_sj_s}-t_{i_2j_1}\cdots t_{i_sj_{s-1}t_{i_1j_s}} $$ which is  in  $I_2(t)$. 

\begin{theorem}(Sturmfels-Villarreal) 
\label{StuVil} 
The set of the polynomials  $f_{I,J}$ where $(I,J)$ is a cycle of $K_{n,m}$ form a
universal  Gr\"obner basis of $I_2(t)$. 
\end{theorem}

In particular we have: 

\begin{corollary} 
\label{StuVil1}
The polynomials $f_{I,J}$ involving only variables of $T(V)$ form a universal Gr\"obner basis of $I_2(t)\cap T(V)$.  
\end{corollary}

Important for us is the following: 

\begin{corollary}  
\label{klaklo}   
The ideal $I_2(t)\cap T(V)$ has a universal Gr\"obner basis whose elements have
$\ZZ^m$-degree bounded above by $(1,1,\dots,1)\in \ZZ^m$. 
\end{corollary}

For a  $\ZZ^m$-graded algebra $E$  we denote by $E_\Delta$    the direct sum of the graded components of
$E$ of  degree $(v,v,\dots,v)\in \ZZ^m$  as $v$ varies in $\ZZ$. Similarly, for a    $\ZZ^m$-graded 
$E$-module  $M$ we denote by  $M_\Delta$   the direct sum of the graded components of
$M$ of  degree $(v,v,\dots,v)\in \ZZ^m$  as $v$ varies in $\ZZ$.  
 Clearly  $E_\Delta$ is a $\ZZ$-graded algebra and $M_\Delta$ is a   $\ZZ$-graded $E_\Delta$-module. Furthermore    $-_\Delta$ is
exact as a functor on the category of $\ZZ^m$-graded $E$-module with maps of degree $0$. 

Now $B(V)_\Delta$ is the $K$-algebra generated by the elements in $y_1V_1\cdots y_mV_m$. Therefore $A(V)$ is (isomorphic to)  the 
algebra of
$B(V)_\Delta$.  

Hence we obtain a presentation 

$$0\to Q \to  T^* \to A(V)\to 0$$
where $Q=(I_2(t)\cap T(V))_\Delta$ and   $T^*=T(V)_\Delta$  is the $K$-algebra generated by the monomials
$t_{1j_1}\dots t_{mj_m}$ with $j_k\in C_k$, that is,   $T^*$ is the  Segre product of the polynomial rings 
$$T_i=K[t_{ij} : j\in C_i].$$
 From \ref{klaklo} we get:

\begin{corollary}
\label{gary} 
The ideal $Q$ is generated by elements of degree $(1,1,\dots,1)$ which form   a Gr\"obner basis with respect to any
term order on  the variables $t_{ij}$. 
\end{corollary}

\begin{proof} Let $g\in Q$ be a  homogeneous element of degree, say, $(a,a,\dots,a)$. Then there exists $h\in I_2(t)\cap T(V)$ of
multidegree
$\leq (1,1,\dots,1)$ such that $\ini(h) | \ini(g)$. Then there exists a monomial $v$ of multidegree
$(1,1,\dots,1)-\deg h$ such that $\ini(h)v |\ini(g)$. It follows that $hv\in Q$ has degree
$(1,1,\dots,1)$  and its initial term divides $\ini(g)$. 
\end{proof} 

In \ref{gary} (and  later on) we consider Gr\"obner bases and initial ideals  of ideals in $K$-subalgebras of polynomial rings.  For the
details on this ``relative" Gr\"obner basis theory the reader can consult,  for instance, 
\cite[Sect.3]{BC} or the \cite[Chap.11]{S}. We may now conclude:

\begin{theorem} 
\label{kosmon}
If the $V_i$ are generated by variables then $A(V)$ is a Koszul  algebra. Moreover
$A(V)$  is a quotient of  the Segre product  $T^*$  by  an ideal generated by   linear  (binomial) forms which are  a Gr\"obner basis. 
\end{theorem}

\begin{proof}    From
\ref{gary} we know that the initial ideal $\ini(Q)$  (with respect to any term order) is an ideal of 
$T^*$ generated by  a subset of the  monomials  generating $T^*$ as a
$K$-algebra.   By work of Herzog, Hibi and Restuccia \cite[2.3]{HHR} we know that Segre products of polynomial rings are
strongly Koszul semigroup rings. Strongly Koszul semigroup rings remain  strongly Koszul after moding out semigroup generators
\cite[2.1]{HHR}. So
$T^*/\ini(Q)$  is strongly  Koszul and in particular Koszul. But then the standard deformation argument shows that $T^*/Q$ is 
Koszul, see
\cite[3.16]{BC} for details. Therefore we can conclude that $A(V)$ is a Koszul algebra. 
\end{proof} 

\begin{remark}
\label{rema1}
  In the proof of above we have shown that a Segre product of polynomial rings modulo a certain ideal of linear forms
is Koszul. One might ask whether linear sections of Segre product of polynomial rings are always Koszul. It is not the case. The
ideal of $2$-minors of the matrix 
$$\left(\begin{array}{cccc}  0 & x & y & z\\ x & y & 0 & t
\end{array}
\right)
$$  defines an algebra which is  a linear section of the Segre product  of polynomial rings of dimension $2$ and $4$ and it is not
Koszul. This  is the algebra number
$69$ in Roos'  list \cite{R}, a well-known gold-mine of examples. 
\end{remark} 

Keeping track of the various steps of the construction above one can describe  the defining equations of
$A(V)$. In details, we set  $C=C_1\times C_2 \times \dots \times C_m$. Consider variables $s_\alpha$ with $\alpha\in C$ and the
polynomial ring
$K[C]=K[s_\alpha : \alpha\in C]$.  Then we get presentations of the Segre product
$T^*$ and of $A(V)$  as  quotients of $K[C]$  of by sending
$s_{(j_1,\dots,j_m)}$ to  $t_{1j_1}\cdots t_{mj_m}$ and to $x_{j_1}\cdots x_{j_m}$ respectively. 

The ring  $T^*$ is the Hibi ring of the distributive lattice $C$ so it is defined by the Hibi relations, namely 

$$s_\alpha s_\beta - s_{\alpha\vee \beta}s_{\alpha\wedge  \beta}  $$

where 

$$\alpha\vee \beta=(\max(\alpha_1,\beta_1), \dots, \max(\alpha_m,\beta_m))$$

and 

$$\alpha\wedge  \beta=(\min(\alpha_1,\beta_1), \dots, \min(\alpha_m,\beta_m)).$$ We have: 

\begin{proposition} 
\label{pseudoW}  The defining ideal of $A(V)$  as a quotient of the polynomial ring $K[C]$ is generated by the   Hibi relations
$s_\alpha s_\beta - s_{\alpha\vee \beta}s_{\alpha\wedge  \beta}$  and by the relations
$$s_\alpha-s_\beta$$ where $\alpha,\beta\in C$ and one is obtained by the other with a non-trivial permutation.
\end{proposition}  

For instance:

\begin{example}
\label{exa1}
Let $n=3$ and $V_1=\langle x_2,x_3 \rangle , V_2=\langle x_1,x_3 \rangle , V_3=\langle x_1,x_2\rangle$. Then
$B(V)$ is the quotient of
$K[t_{12}, t_{13},t_{21}, t_{23}, t_{31}, t_{32}]$ by the polynomial
$t_{12}t_{23}t_{31}-t_{13}t_{21}t_{32}$ and  then
$A(V)$ is the quotient of 
$K[s_{ijk} : (i,j,k)\in \{2,3\}\times \{1,3\}\times \{1,2\}]$  by the Hibi-relations

$$\begin{array}{ll}
  s_{312}s_{331}-s_{311}s_{332}, &
  s_{212}s_{311}-s_{211}s_{312},\\
  s_{212}s_{231}-s_{211}s_{232}, &
  s_{212}s_{331}-s_{211}s_{332},\\
  s_{231}s_{311}-s_{211}s_{331}, &
  s_{231}s_{312}-s_{211}s_{332},\\
  s_{232}s_{311}-s_{211}s_{332}, &
  s_{232}s_{312}-s_{212}s_{332},\\
  s_{232}s_{331}-s_{231}s_{332}
\end{array}
$$

and by the linear relation

 $$s_{231}-s_{312}$$
\end{example}

\begin{remark}
\label{rema2}
  It is not clear whether the defining ideal of $A(V)$ as a quotient of $K[C]$ has a Gr\"obner basis of quadrics. 
 The Hibi relations form a Gr\"obner basis  with respect  to any  revlex linear extension   of the
partial order on $C$.  There are examples where the Hibi relations together with the linear relations  defining $A(V)$ are not a
Gr\"obner basis with respect to such revlex linear extensions. 
\end{remark} 

\begin{remark}
\label{rema3}
In a special case it turns out that both   $B(V)$ and $A(V)$ are defined by   Gr\"obner bases of quadrics as quotient
of polynomial rings.  For a nested chain of vector spaces of linear forms 
$V_1\supseteq V_2\supseteq \dots \supseteq V_m$,  we can fix a basis $x_1,x_2,\dots,x_n$ of $R_1$ such that $V_i$ is
generated by
$x_1,\dots,x_{d_i}$. Here   $d_1\geq d_2\geq \dots\geq d_m$.  It follows that
$B(V)$ corresponds to a  one-sided ladder determinantal ring, the ladder being the set of points  $(i,j)$ with $1\leq i\leq m$  and
$1\leq j\leq d_i$. Furthermore, $A(V)$ coincides with the algebra associated with the principal Borel subset generated by the
monomial $\Pi_i x_{d_i}$.  A Gr\"obner basis of quadrics  for
$B(V)$ is  described in \cite{HT} and  a Gr\"obner basis of quadrics  for $A(V)$ is  described in \cite{D}. 
\end{remark} 

In general, however, the algebra $B(V)$ is not defined by quadrics as the Example \ref{exa1} shows. White's conjecture \cite{Wh}
predicts the structure of the defining equations of the base ring of a (poly)matroid: they should be quadrics representing the basic
symmetric exchange relations of the polymatroid. Our result  above \ref{pseudoW} does not prove White's  conjecture in this
precise form.

\section{Conjectures}  
\label{conj}
The constructions and arguments of the previous section suggest a general strategy to investigate the Koszul
property of
$A(V)$ for general (i.e. non-monomial)  $V_i$. We outline in this section the strategy which leads us to Conjecture \ref{strocon}. 
Let $V=V_1,\dots,V_m$ be a collection of subspaces of $R_1$ and let $y_1,\dots, y_m$ be new variables.  Set 
$d_i=\dim V_i$,  and set 

$$S=K[y_ix_j :  i=1,\dots,m, j=1,\dots,n]$$

$$B(V)=K[y_1V_1,\dots,y_mV_m].$$

and
 
$$T=K[t_{ij} :  i=1,\dots,m, j=1,\dots,n].$$

Again $B(V)$ is a
$K$-subalgebra of  $S$. We give degree
$e_i\in {\bf Z}^m$ to $y_ix_j$ and to $t_{ij}$ so that  
$S$, $T$  and as well $B(V)$ are ${\bf Z}^m$-graded.  We present $S$ as a quotient of $T$ by sending 
$t_{ij}$ to $y_ix_j$. The kernel of such presentation is the ideal $I_2(t)$ generated by the $2$-minors of the $m\times n$ matrix
$t=(t_{ij})$.  As we have   seen in the previous section $A(V)$ is the diagonal algebra $B(V)_\Delta$.

We want to get the presentations of $B(V)$ and $A(V)$ by elimination from that of $S$. To that end  we do the following:  Let $f_{ij}$,
$j=1,\dots,d_i$,   be a basis of $V_i$ and complete it to a basis of $R_1$ with elements  $f_{ij}$, $j=d_i+1,\dots,n$. Denote by $f_i$
the row  vector $(f_{ij})$ and by $x$  the row vector of the
$x_i$'s. Let  $A_i$ be the $n\times n$ matrix with entries in $K$ with $x=f_i A_i$.  Then 
$S=K[y_if_{ij} : i=1,\dots,m,\  j=1,\dots,n]$ and
$B(V)=K[y_if_{ij} : i=1,\dots,m,\  j=1,\dots,d_i]$. Set $T(V)=K[t_{ij} : 1\leq i\leq m, 1\leq j\leq d_i]$. We have presentations: 

$$\begin{array}{lll} 
\phi:  T     \to S &      \mbox{  with }  t_{ij}\to y_if_{ij}  &  \mbox{ for all } i,j \\  \\ 
\phi': T(V)\to B(V) &      \mbox{  with }  t_{ij}\to y_if_{ij}  & \mbox{ for all } i \mbox{ and } 1\leq j\leq d_i
\end{array}
$$

By construction, the kernel of  $\phi$ is the ideal of $2$-minors $I_2(L)$ of the matrix 
$L=(L_{ij})$ where the  row vector $(L_{ij} : j=1,\dots,n)$ is  given by $(t_{i1},\dots, t_{in})A_i$.  Clearly, $\Ker \phi'=I_2(L)\cap T(V)$. 
As explained in the previous section, by applying the diagonal functor  we obtain  a presentation: 

$$A(V)\simeq T^*/Q$$
where $T^*$ is the Segre product of  the $T_i$'s,  $T_i=K[t_{ij}  : j=1,\dots,d_i]$,  and $Q=(I_2(L)\cap T(V))_\Delta$.

\begin{remark}
\label{rema4}
 One can easily check that the arguments of Section \ref{mono}, in particular those of \ref{gary} and \ref{kosmon},   work    and  can
be used to show that $A(V)$ is Koszul provided one knows that
$I_2(L)\cap T(V)$ has an initial ideal generated in degree $\leq (1,1,\dots,1)\in \ZZ^m$.  On the other hand,  $I_2(L)\cap T(V)$ has the
desired initial ideal provided  $I_2(L)$ has an initial ideal 
 generated in degree $\leq (1,1,\dots,1)\in \ZZ^m$ with respect to the appropriate elimination order. 
\end{remark}

We are  led by \ref{rema4} to analyze initial ideals of ideals of $2$-minors of matrices as $L$. To our great surprise, the
experiments   support the  Conjecture \ref{strocon}.   What we really need is a weak for of \ref{strocon},
namely: 

\begin{conjecture}
\label{weakcon} 
  Let $L=(L_{ij})$ be a $m\times n$ matrix with $L_{ij}=\sum_{k=1}^n a_{ijk}t_{ik}$
and
$a_{ijk}\in K$ for all $i,j,k$. Assume that for every $i$ the forms
$L_{i1},\dots,L_{in}$ are linearly independent.  Then  any lexicographic initial ideal of
$I_2(L)$ is generated in degree $\leq (1,1,\dots,1)$.  
\end{conjecture} 

 If conjecture \ref{weakcon} holds then  from  the discussion above follows that     for every $V_1,\dots, V_m$
the algebra $A(V)$ is Koszul and defined by a Gr\"obner basis of  linear forms  as a quotient of the Segre product $T^*$. 

The next section is devoted to prove Conjecture \ref{strocon} in the generic case.

\section{The  generic case}
\label{gene}

We  consider  now the case of generic spaces $V_1,\dots,V_m$. What we prove is the following: 

\begin{theorem}
\label{main} 
If the matrix $L$ is generic, that is, every entry $L_{ij}=\sum_{k=1}^n a_{ijk}t_{ik}$  is a generic linear
combination of the variable $t_{i1}, \dots,t_{in}$, then 
\ref{strocon} holds. 
\end{theorem}

The key lemma is:  

\begin{lemma}
\label{key}
 Let  $V_1,\dots,V_m$  be  subspaces of $R_1$.   If   $\sum_{i=1}^m  \dim V_i \geq n+m$  then $\dim
\prod_{i=1}^m V_i < \prod_{i=1}^m  \dim V_i$, i.e.   there is a non-trivial linear relation  among the  generators of the
product  $\prod_{i=1}^m V_i$ obtained by multiplying $K$-bases of the $V_i$.
\end{lemma} 

\begin{proof}  By induction on $n$ and $m$. If one of the $V_i$ is principal then we can simply skip it. The case $m=2$ is easy: the
assumption is equivalent to $\dim(V_1\cap V_2)\geq 2$ and for $f,g\in V_1\cap V_2$ we get the non-trivial relation
$fg-gf=0$. For  $m>2$,  if $\dim(V_i\cap V_j)\geq 2$ for some $i\neq j$ then the non-trivial relation above gives a non-trivial
relation also for
$V_1\cdots V_m$. Therefore we may assume that   $\dim(V_i\cap V_j)< 2$, and, since none of the $V_i$ is principal, also none of the
$V_i$ is $R_1$. The case $n=2$ follows and  to prove the assertion in the general case we may assume that  $1< d_i <n$ for all $i$.
Further we may assume also that the
$V_i$ are generic, the dimension of $V_1\cdots V_m$ for special $V_i$ can be only smaller. By the genericity of the $V_i$ we may find
$K$-bases $f_{ij}$ of
$V_i$ so that any set of $n$ elements in the  set  $\{ f_{ij} : i=1,\dots,m, \hbox{ and }j=1,\dots, d_i\}$ is a basis of $R_1$. Now let $x$
be a general linear form (it suffices that $x$ is not contained in any sum of the $V_i$ which is a proper subspace of $R_1$).  Since
$x\not\in V_i$ we  have that $\dim V_i+(x)/(x)=d_i$, so by induction on $n$ we may find a non-trivial relation among the  
generators of
$V_1\cdots V_m$ modulo
$x$. In other  words there exists a relation of the form  
$$\sum \lambda_{\alpha}f_{1\alpha_1}\cdots f_{m\alpha_m}=xh$$  where $\lambda_{\alpha}\in K$, the sum is extended to all the
$\alpha$ in 
$\prod_{i=1}^m \{1,\dots,d_i\}$ and at least one of the $\lambda_{\alpha}$ is non-zero.  We may assume $\lambda_{\alpha}\neq 0$
 for 
$\alpha=(1,1,\dots,1)$. By the above relation we have that $xh\in \prod_{i=1}^m  V_i$ and hence 
$xh\in \prod_{i\neq j} V_i$ for all $j$. But form \ref{primdec}  we see immediately that $x$ acts as a non-zero divisor  in  degree
$m-1$  and higher  on  the  ideal generated by 
$\prod_{i\neq j} V_i$.  It follows that $h\in
\prod_{i\neq j} V_i$ for all $j$.  By the choice of the $f_{ij}$ and since $\sum_{i=1}^m  d_i\geq n+m$  we may write $x$ as a linear
combination   of the $f_{ij}$ with $i=1,\dots,m, $  and  $1<j\leq d_1$. It follows that
$xh$ can be written as a linear combination of  the $f_{1\alpha_1}\cdots f_{m\alpha_m}$ with $\alpha\neq (1,1,\dots,1)$. Hence we
obtain a relation  
$$\sum \lambda_{\alpha}^\prime f_{1\alpha_1}\cdots f_{m\alpha_m}=0$$ with  $\lambda_{\alpha}^\prime=\lambda_{\alpha}\neq
0$ for $\alpha= (1,1,\dots,1)$.
\end{proof} 

Now we are ready to prove:

\begin{proof} of \ref{main} Set $I=I_2(L)$.  Let $<$ be a term order on the $t_{ij}$. After a change of name of the variables in the
$i$-th row of
$L$ if needed, we may assume that  $t_{ij+1}>t_{ij}$ for all
$j=1,\dots,n-1$ and for all $i=1,\dots,m$.  Let $J$ be the ideal generated by   the monomials  

$$t_{i_1j_1}\cdots t_{i_kj_k}$$

satisfying conditions: 

$$ (*) \left\{\begin{array}{l} 1\leq i_1<\dots<i_k\leq m,\\ 1\leq j_1,\dots, j_k\leq n,\\ j_1+\dots+j_k\geq n+k. 
\end{array}\right.
$$ 

We will show that the initial ideal  of $I$ with respect to $<$ is equal to  $J$. From this the assertion follows  immediately.   
It is a simple exercise on primary decompositions that the equality 
$J=\ini(I)$  follows from   three  facts:\\   (1)  $J\subseteq \ini(I)$, \\    (2)  $J$ and $I$ have the same codimension and degree,\\ (3) 
$J$ is unmixed.\\
 
For (1)  we  have to show that for each  pair of sequences of integers satisfying conditions (*) 
 the monomial  $t_{i_1j_1}\cdots t_{i_kj_k}$ is in $\ini(I)$. As $L$ is generic, the initial ideal
$\ini(I)$ is the multigraded generic initial ideal  of
$I$ with respect to $>$. Hence $\ini(I)$  is Borel fixed is the multigraded sense, see
\cite{ACD}. In characteristic $0$ this means that if  a monomial $M$ is in $\ini(I)$ and
$t_{ij}| M$ then  $t_{ik}M/t_{ij}$ is in $\ini(I)$ as well for all the $k>j$.  In arbitrary characteristic the same assertion is also  true 
as long as $M$ is square-free.  It follows that,  (no matter what the characteristic is),  it suffices to show that there exists an  $f$ in
$I$ such that
$\ini(f)=t_{i_1p_1}\cdots t_{i_kp_k}$  and
$p_1\leq j_1,\dots, p_k\leq j_k$.  To this end, consider the linear forms 
$f_{ij}$ defined (implicitly)  by the relation $x_k=\sum_{k=1}^n f_{ij} a_{ikj}$ for all $j$.  By the construction of Section \ref{conj}
we see that $I$ is the kernel of the map $\phi$.  Now for  $s=1,\dots,k$ consider the subspace
$W_{i_s}$  generated by  the $f_{i_sj}$  with $j\leq j_s$. Since,   by assumption 
$\sum_{s=1}^k \dim W_{i_s}=\sum_{s=1}^k  j_{s}\geq n+k$, by Lemma \ref{key}  we have that there exists a non-trivial relation
among the   generators of  the product 
$W_{i_1}\cdots W_{i_k}$. This implies that the $I$ contains a non-zero polynomial
$f$ supported on the set of the monomials $t_{i_1p_1}\cdots t_{i_kp_k}$  and
$p_1\leq j_1,\dots, p_k\leq j_k$. Take $\ini(f)$ to get what  we want.

As for the step (2) and (3),  the ideal of $I$ is a generic determinantal ideal  and its numerical invariants  are well-known:   its
codimension  is
$(m-1)(n-1)$ and the its degree is
$\binom{m+n-2}{m-1}$.  Knowing the generators of
$J$ we can describe the facets of the associated simplicial complex $\Delta(J)$. Then we can read from the descriptions of the
facets  the codimension, the degree of $J$ and  check that it is unmixed. The facets of $\Delta(J)$ have the following description: for
each $p=(p_1,\dots,p_m) \in \{1,\dots,n\}^m$ with $p_1+\dots,p_m=n+m-1$ we let 

$$F_p=\{ t_{ij} : i=1,\dots, m \mbox{ and }  1\leq j\leq p_i\}$$

It is easy to   check that any such $F_p$ is a facet of $\Delta(J)$. On the other hand if $F$ is a face of
$\Delta(J)$ let $a(F)=\{ i : \exists\  j \mbox{ with } t_{ij}\in F\}$ and $j_i=\max\{ j : t_{ij}\in F\}$ if $i\in a(F)$. Then set
$q=(q_1,\dots,q_m)$ with
$q_i=j_i$ if $a\in a(F)$ and $q_i=1$ otherwise.  Note that  
$$q_1+\dots+q_m=\sum_{i\in a(F)} j_i +m-|a(F)|$$
 and that 
$$ \sum_{i\in a(F)} j_i <n+|a(F)|$$
 since $\{ t_{ij_i} : i \in a(F)\}\subset F\in
\Delta(J)$. It follows that  $q_1+\dots+q_m<n+m$. So, increasing the $q_i$'s if needed,  we may take
$p=(p_1,\dots,p_m)\in 
\{1,\dots,n\}^m$ with $p_1+\dots,p_m=n+m-1$ and $q_i\leq p_i$. It follows that $F\subseteq F_p$. 

From the description above we see that the cardinality of each $F_p$ is $n+m-1$. It follows that $J$ is unmixed of codimension
$(m-1)(n-1)$. 
 The degree $J$ is the number of facets of $\Delta(J)$, that is the number of  $p=(p_1,\dots,p_m) \in
\{1,\dots,n\}^m$ with $p_1+\dots+p_m=n+m-1$. Setting $q_i=p_i-1$, we see that the number of facets of
$\Delta(J)$ is the number of $q=(q_1,\dots, q_m)\in \{0,\dots,n-1\}^m$ with $q_1+\dots+q_m=n-1$, that is, the number of monomials
of degree $n-1$ in $m$ variables. This number is $\binom{m+n-2}{m-1}$.   
 We have checked that (2) and (3) hold.  The proof of the theorem is now  complete. 
\end{proof}

Let us single out the following corollary of the proof of of \ref{main}:

\begin{corollary}
\label{equat} 
With the notations of the proof of \ref{main} we have:
\begin{itemize}
\item[(a)]  If $i_1<\dots<i_k$ then a monomial $t_{i_1j_1}\cdots t_{i_kj_k}$ is in $J$ iff  $j_1+\dots+j_k\geq n+k$. 
\item[(b)] For every monomial $M=t_{i_1j_1}\cdots t_{i_kj_k}\in J$  with  $i_1<\dots<i_k$    there exists a polynomial $f_M\in I$ of
the form 
$$f_M=M+\sum_v \lambda_v t_{i_1v_1}\cdots t_{i_kv_k}$$ where $\lambda_v\in K$,  $v\in \Pi_{h=1}^k \{1,2,\dots,j_h\}$, and 
$t_{i_1v_1}\cdots t_{i_kv_k}\not\in J$.
\item[(c)]   The set of the polynomials $f_M$ is  a  Gr\"obner basis of $I$ with respect to any term order $<$ on $K[t_{ij}]$ satisfying 
$t_{ij+1}>t_{ij}$ for all $j=1,\dots,n-1$ and for all $i=1,\dots,m$.
\end{itemize}
\end{corollary} 

\begin{proof} (a) follows form the definition of $J$. For (b) we argue as follows. Let  $<$  be a term order on $K[t_{ij}]$ satisfying 
$t_{ij+1}>t_{ij}$ for all
$j=1,\dots,n-1$ and for all $i=1,\dots,m$. We have seen in the proof of \ref{main} that  
$J=\ini_<(I)$.  Considering the reduced expression, we have that for every  monomial $M=t_{i_1j_1}\cdots t_{i_kj_k}\in J$  there
exists a polynomial
$f_M$ in $I$ with  initial  term $M$ and   all the others terms not in $J$. Suppose that one the non-leading terms of
$f_M$, say $N=t_{i_1v_1}\cdots t_{i_kv_k}$, does not satisfies the  condition $v_h\leq j_h$ for all
$h=1,\dots,k$.  So there exists  an $h$ in $\{1,2,\dots,k\}$, say $h_1$,  such that
$v_{h_1}>j_{h_1}$.  We claim that  there exists a term order $<_1$ such that 
$t_{ij+1}>_1t_{ij}$ for all $i,j$ and such that  $N>_1M$. Then it follows that the initial term of $f_M$ with respect to $<_1$ is not $M$
and hence in must be a monomial  not in $J$. This contradicts the fact, proved in \ref{main} that $\ini_{<_1}(I)=J$.  It remains to
prove the existence of a term order
$<_1$ as above.  To this end  it  is suffices to find weights $w_{ij}\in \NN$ such that $w_{ij}<w_{ij+1}$ for all
$i,j$ and 
$w(M)<w(N)$, that is 

$$w_{i_1,j_1}+\dots + w_{i_kj_k} < w_{i_1v_1}+\dots + w_{i_kv_k}$$

Just take $w_{ij}=j$ if $i\neq i_{h_1}$  of if $i= i_{h_1}$ and $j<v_{h_1}$ and $w_{ij}=a+j$ otherwise  with $a$ large enough. Finally
(c) is a direct consequence of  (b). 
\end{proof}

As explained in Section  \ref{conj}  from \ref{main} follows that $A(V)$ is Koszul for generic $V$. To get more  precise information
about the structure of $A(V)$ we  analyze in  details  the defining equations of of $B(V)$ and $A(V)$.   
To this end we recall the definition of homogeneous ASL on a poset.

Let $(H,>)$ be a finite poset and denote by $K[H]$ the polynomial  ring whose variables are  the elements of $H$.  
Let $J_H$ be the monomial ideal of $K[H]$ generated by $xy$ with $x,y\in H$ such that $x$ and $y$ are incomparable in $H$. 

\begin{definition}
\label{ASL}
 Let $A=K[H]/I$ where  $I$  is a homogeneous ideal (with respect to the usual grading). One says that $A$
is a homogeneous  ASL on $H$ if 
\begin{itemize}
\item[(ASL1)] The (residue classes of the) monomials not in $J_H$ are linearly independent in $A$. 
\item[(ASL2)]  For every $x,y\in H$ such that $x$ and $y$ are incomparable the ideal $I$ contains a polynomial of the form 
$$xy-\sum \lambda zt$$
 with  $\lambda \in K$,  $z, t\in H$, $z\leq t$, $ z<x$   and  $z<y $. 
\end{itemize} 
\end{definition}

 A  linear extension of  the poset $(H,<)$  is a total order $<_1$ on $H$ such that $x<_1 y$ if $x<y$. 
A revlex term order $\tau$ on $K[H]$ is said to be a revlex linear extension of $<$ if $\tau$ induces on $H$ a linear extension of $<$.  
For obvious reasons, if $A=K[H]/I$  is a homogeneous  ASL 
on $H$ and $\tau$ is a  revlex linear extension of $<$ then the polynomials in (ASL2) form a Gr\"obner basis of $I$   and 
$\ini_{\tau}(I)=J_H$.  In a  sense the converse in also   true:

\begin{lemma} 
\label{ASL-GB}
Let $A=K[H]/I$ where  $I$  is a homogeneous ideal. Assume that for
every revlex linear extension $\tau$ of $<$ one has $\ini_{\tau}(I)=J_H$. Then $A$ is an ASL on $H$. 
\end{lemma}  

\begin{proof} Let $\tau$ be a  revlex linear extension of $<$. Since $\ini_{\tau}(I)=J_H$  the monomials  not in $J_H$ form a
$K$-basis of $A$, hence ASL1 is satisfied.  Let $x,y\in H$ be incomparable elements. Then $xy\in \ini_{\tau}(I)$ and hence there exists
$F\in I$ with
$\ini_{\tau}(F)=xy$. We can take $F$ reduced  in the sense that $xy$ is the only term in $F$ belonging to $J_H$. It follows that  $F$
have the form 
$$xy-\sum \lambda zt $$
  with  $\lambda \in K$, $z, t\in H$ and $z\leq t$. Assume, by contradiction that this polynomial does not satisfy the conditions
required in ASL2.  Then there exist a  non-leading  term $z_1t_1$ appearing in $F$  such that  $z_1\not< x$ or  $z_1\not< y$. Say  
$z_1\not< x$.
 It is easy to see that one can find a linear extension $<_1$ of $<$ such that $x<_1 z_1$.  Denote by $\sigma$ the revlex term order 
associated with $<_1$. 
Then   $xy$ is smaller than $z_1t_1$ with respect to $\sigma$ and hence
$\ini_{\sigma}(F)$ is a term not in
$J_H$, contradicting the assumption.    
\end{proof}

For a given sequence of positive integers $d=d_1,\dots,d_m$ we set 

$$H(d)=\{1,\dots,d_1\}\times \dots \times \{1,\dots,d_m\}$$

and note that $H(d)$ is a sublattice of $\NN^m$ with respect to the natural partial order    
$\alpha\leq \beta$ iff $\alpha_i\leq \beta_i$ for all $i$. The rank $\rk \alpha$  of an element
$\alpha=(\alpha_i) \in H(d)$ is $\alpha_1+\dots+\alpha_m-m$. 

Set 

$$H_n(d)=\{ \alpha \in H(d) : \rk \alpha<n\}$$ 

With the notation of Section \ref{conj}  we have a presentation  $\phi':T(V)\to B(V)$  where 
$T(V)=K[t_{ij} :  i=1,\dots,m, \ \ j=1,\dots, d_i]$. 
 As a corollary of \ref{main}, by elimination we obtain  a description  $\Ker \phi'$: 
 
\begin{corollary} 
\label{dama1}
 Let $V_1,\dots,V_m$ be generic spaces of dimension $d_1,\dots,d_m$ and let $f_{ij}$ with $j=1,\dots,d_i$
be generic generators of
$V_i$.  Let $<$ be a term order such that $t_{ij}<t_{ij+1}$. Then the ideal 
$\Ker\phi'$ has a Gr\"obner basis whose elements are the polynomials $f_M$ of \ref{equat} with
$M=t_{i_1j_1}\cdots t_{i_kj_k}$  with  $i_1<\dots<i_k$,  $1\leq j_h\leq d_{i_h}$ and $j_1+\dots+j_k\geq n+k$.  \end{corollary} 

Set
$T_i=K[t_{ij} : 1\leq j\leq d_i]$ and denote by $T^*$  the Segre product
$T_1*\dots *T_m$.  Consider variables $s_\alpha$ with $\alpha\in H(d)$ and  the polynomial  ring $K[s_\alpha : \alpha\in H(d)]$.   
  For each
$\alpha\in H(d)$  set $t_\alpha=t_{1\alpha_1}\cdots t_{m\alpha_m}$.

We get a presentation   $K[s_\alpha : \alpha\in H(d)]\to T^*$ by sending $s_\alpha$ to
$t_\alpha$ whose kernel is generated by the Hibi relations:

$$s_\alpha s_\beta-s_{\alpha\vee \beta}s_{\alpha\wedge \beta}.$$

Adopting the notation of Section \ref{conj} we get a presentation $A(V)=T^*/Q$. To describe the generators of $Q$ we do the
following. For every
$\alpha \in H(d)\setminus H_n(d)$  consider the polynomial $f_M$ of \ref{equat} associated with the monomial
$M=t_{1\alpha_1}\cdots t_{m\alpha_m}$.  Set $L_\alpha=f_M$.  So for all $\alpha\in H(d)\setminus H_n(d)$ we have 
$$L_\alpha=t_\alpha-\sum_{\beta<\alpha} \lambda_{\alpha\beta} t_\beta  \quad \mbox{ with }
\lambda_{\alpha\beta}\in K$$  and  the arguments of \ref{gary} show that the $L_\alpha$'s  form a Gr\"obner basis of $Q$ for any
term order such that
$t_{ij}>t_{ij-1}$ for all $i,j$. It follows that

$$\ini(Q)=(t_{1\alpha_1}\cdots t_{m \alpha_m} :  \alpha \in H(d)\setminus H_n(d)) $$

for any term order such that $t_{ij}>t_{ij-1}$ for all $i,j$. Then $T^*/\ini(Q)$ is defined as a quotient of $K[s_\alpha : \alpha\in H(d)]$
by: 
\medskip

\nt (1) the Hibi relations $s_\alpha s_\beta-s_{\alpha\vee \beta}s_{\alpha\wedge \beta}$ with
$\alpha,\beta\in H(d)$ incomparable. 
\medskip 

 \nt (2) $s_\alpha$ with $\alpha\in H(d)\setminus H_n(d)$. 
\medskip 

It is easy to see that the elements of type  (1) and (2) form a Gr\"obner basis for any revlex linear extension of the partial order on
$H(d)$. Hence a $K$-basis of  $T^*/\ini(Q)$ is given by  the monomials not in $J_{H_n(d)}+(  H(d)\setminus H_n(d))$. 
This in turns implies that the Hibi relations and the relation
$L_\alpha$ form a Gr\"obner basis with respect to any revlex linear extension of the partial order on $H(d)$ of the defining ideal of
$A(V)$ as a quotient of
$K[s_\alpha : \alpha\in H(d)]$ by the map sending
$s_\alpha$ to $f_{1\alpha_1}\cdots f_{m\alpha_m}$.  Summing up, we have:

\begin{theorem}
\label{ASLA}
 Let $V_1,\dots,V_m$ be generic spaces of dimension $d_1,\dots,d_m$  and take generic generators $f_{ij}$ of $V_i$. Then: 
\begin{itemize} 
\item[(1)]  We have a surjective
$K$-algebra  homomorphism  $F:K[s_\alpha : \alpha\in H_n(d)]\to A(V)$ sending the variable $s_\alpha$ to
$f_{1\alpha_1}\cdots f_{m\alpha_m}$.
\item[(2)] $\Ker F$ is generated by two types of polynomials: 
\begin{itemize} 
\item[(a)]  $$s_\alpha s_\beta-s_{\alpha\vee \beta}s_{\alpha\wedge\beta}$$  if $\alpha,\beta\in H_n(d)$ are incomparable and
$\alpha\vee \beta \in H_n(d)$. 
\item[(b)]  $$s_\alpha s_\beta- \sum \lambda_\gamma  s_\gamma s_{\alpha\wedge \beta}$$  if $\alpha,\beta\in H_n(d)$ are
incomparable and
$\alpha\vee \beta \not\in H_n(d)$  and  the sum is extended to the $\gamma\in H_n(d)$ with $\gamma\leq \alpha\vee
\beta$ and
$\lambda_\gamma\in K$ (and depends also on $\alpha$ and $\beta$).
\end{itemize}  
\item[(3)] The polynomials of type (a) and (b) form a Gr\"obner basis of $\Ker F$ with respect to any revlex  linear extension of 
 the partial order of  $H_n(d)$. 
\item[(4)]  $A(V)$ is a homogeneous ASL on the  poset $H_n(d)$. 
\item[(5)] $A(V)$ is normal, Cohen-Macaulay and  Koszul. 
\item[(6)]  $A(V)$ is defined,   as a quotient of the Segre product $T^*$, by a Gr\"obner basis of linear forms. 
\item[(7)]  The Krull dimension  of $A(V)$ is $\min\{n, \dim T^*=1-m+\sum_{i=1}^m  d_i \}$ and its degree is the number of
maximal chains in $H_n(d)$. 
\end{itemize} 
\end{theorem}

\begin{proof} (1), (2), (3) and (6) follows immediately from the discussion above and  (4) follows \ref{ASL-GB} and (3).   As for (5),
normality is proved in \ref{normal}, Koszulness follows from the general argument of Section
\ref{conj} and also from  (3). The Cohen-Macaulay property and (7) follows from (4) by applying
\cite[Chap.5]{BV} since $H_n(d)$ is a wonderful poset.   
\end{proof}

As a corollary we obtain: 

\begin{corollary}
\label{noamdg}
  For every $m$ and $n$, the Veronese subring $R^{(m)}$ of $R=K[x_1,\dots,x_n]$ is an ASL on the poset 
$H_n(d)$ where $d=n,n,\dots, n$ ($m$-times). 
\end{corollary} 

\begin{remark} 
\label{larema}
 The realization of the $m$-th Veronese subring of a polynomial ring in $n$ variables as a homogeneous ASL has been
done before for $n=2$ and any $m$ in \cite{W1},  for  $n=m=3$ in  \cite{HW1} and in two different ways, and for $n=m=4$ in
\cite{W2}.
\end{remark}

 An  interesting consequence of \ref{dama1} is: 

\begin{corollary}
\label{lacora}  Let $V_1,\dots,V_m$ be subspaces  of $R_1$ of dimension $d_1,d_2,\dots,d_m$ then:
\begin{itemize} 
\item[(a)]  $\dim \Pi_{i=1}^m V_i\leq |H_n(d)|$.  
\item[(b)] if the $V_i$ are generic then  $\dim \Pi_{i=1}^m V_i=|H_n(d)|$.
\item[(c)] if the $V_i$ are generic  and  if $f_{ij}$ with
$j=1,\dots,d_i$  are generic generators of $V_i$  then the set $\{ f_{1j_1}\cdots f_{mj_m} :    (j_1,\dots, j_m)\in H_n(d)\}$   is a
$K$-basis of
$\Pi_{i=1}^m V_i$.  
\item[(d)] if the $V_i$ are generic then:  $\dim \Pi_{i=1}^m V_i=\Pi_{i=1}^m \dim V_i$ iff $\sum \dim V_i<m+n$. 
\end{itemize} 
\end{corollary} 

\begin{proof}  Obviously  (b) implies  (a) and also  (c) implies   (b) and (d). So we have only to prove (c). By definition, the product 
$\Pi_{i=1}^m V_i$ is the component of degree $(1,1,\dots,1)$ of the algebra $B(V)$.  Then the conclusion follows from \ref{dama1}. 
\end{proof}  

\begin{example}
\label{exala} Take $n=3$ and $d_1=d_2=d_3=2$ and generic spaces $V_i$ of dimension $d_i$. Note that,  up to a choice of
coordinates,  we are in the situation of Example \ref{exa1} and so the structure of $A(V)$ has been already identified. But
to describe the  the ASL structure of $A(V)$ we  have take generic coordinates for the $V_i$, say $V_i=\langle f_{i1}, f_{i,2}\rangle$.  
In this case
$H_n(d)$ is the cube $\{1,2\}^3$ without the point $(2,2,2)$. We have a  relation 
$$f_{12}f_{22}f_{32}= \sum_{\alpha\in H_n(d)}  \lambda_\alpha  f_{1\alpha_1}f_{2\alpha_2}f_{3\alpha_3}.$$ 
Set $L=\sum_{\alpha\in H_n(d)} 
\lambda_\alpha s_\alpha$.  Then the defining equations of $A(V)$ as a quotient of $K[s_\alpha : \alpha \in H_n(d)]$ are: 

$$\begin{array}{lll} s_{112}s_{221} -s_{111}L, \ \ & s_{121}s_{212} -s_{111}L,  &  \ \ s_{211}s_{122}-s_{111}L, 
\\ s_{121}s_{211} - s_{111}s_{221},  &   s_{112}s_{211} - s_{111}s_{212}, &   s_{112}s_{121} - s_{111}s_{122}, \\
 s_{212}s_{221} - s_{211}L, &  s_{122}s_{221} - s_{121}L,  & s_{122}s_{212} -s_{112}L  
\end{array} 
$$

\end{example} 

\begin{remark}
\label{reB}
 With an argument similar to that of \ref{normal} one can prove that the algebra $B(V)$ is normal for any
$V=V_1,\dots,V_m$. Furthermore, in the monomial and in the generic  case one can  prove that $B(V)$ is Cohen-Macaulay. In the
monomial case the Cohen-Macaulayness  is a consequence of the normality. In  generic case it follows form the fact that, by
\ref{main}, we can describe an  initial ideal of its defining ideal and such initial ideal turns out to be associated with a shellable
simplicial complex.
\end{remark}


\begin{thebibliography}{99}

\bibitem{ACD} A.Aramova,  K.Crona, E.De Negri  {\em   Bigeneric initial ideals, diagonal subalgebras and bigraded Hilbert functions}.  J.
Pure Appl. Algebra 150 (2000),no. 3, 215--235.


\bibitem{BC} W.Bruns,  A.Conca,  {\em  Gr\"obner bases and determinantal ideals}. Commutative algebra, singularities and computer
algebra (Sinaia, 2002), 9--66, NATO Sci. Ser. II Math. Phys. Chem., 115, Kluwer Acad. Publ., Dordrecht, 2003. 

\bibitem{BH}	W. Bruns and J. Herzog, {\em  Cohen-Macaulay Rings},  Cambridge University Press, Cambridge, 1996. 

\bibitem{BV} W.Bruns,  U.Vetter, { \em  Determinantal rings}. Lecture Notes in Mathematics, 1327. Springer-Verlag, Berlin, 1988.  

\bibitem{Cocoa}  CoCoA  Team, {\em    CoCoA: a system for doing 
     Computations in Commutative Algebra},   Available at http://cocoa.dima.unige.it 


\bibitem{CH} A.Conca, J.Herzog {\em Castelnuovo-Mumford regularity of products of ideals}. Collect. Math. 54 (2003), no. 2,
137--152. 

\bibitem{DEP} C.De Concini, D.Eisenbud, C.Procesi, {  \em Hodge algebras},   Ast\'erisque, 91.  Soci\'et\'e Math\'ematique de France,
Paris, 1982. 87 pp.
 
\bibitem{D}  E.De Negri,  {\em Toric rings generated by special stable sets of monomials}. 
 Math. Nachr. 203 (1999), 31--45.
 

\bibitem{Ed}  J. Edmonds, { \em Submodular functions, matroids, and certain polyhedra,}  in Combinatorial Structures and Their
Applications, R. Guy, H. Hanani, N. Sauer, and J. Schonheim (Eds.), Gordon and Breach, New York, 1970, pp. 69--87.

\bibitem{Ei} D.Eisenbud,   {\em Introduction to algebras with straightening laws.}  Ring theory and algebra, III (Proc. Third Conf.,
Univ. Oklahoma, Norman, Okla., 1979), pp. 243--268,  Lecture Notes in Pure and Appl. Math., 55,  Dekker, New York, 1980. 


\bibitem{HH} J.Herzog and T.Hibi, { \em Discrete polymatroids}. J. Algebraic Combin. 16 (2002), no. 3, 239--268.

\bibitem{HHR} J.Herzog and T.Hibi, G.Restuccia {\em  Strongly Koszul algebras}.   Math. Scand. 86 (2000), no. 2, 161--178.

\bibitem{HT} J.Herzog and N.V.Trung, { \em Gr\"obner bases and multiplicity of determinantal and Pfaffian ideals}. Adv. Math. 96
(1992), no. 1, 1--37.  

\bibitem{H1}   T.Hibi,  {\em Distributive lattices, affine semigroup rings and algebras with straightening laws}.
 Commutative algebra and combinatorics (Kyoto, 1985), 93--109,Adv. Stud. Pure Math., 11, North-Holland, Amsterdam, 1987.  
 
\bibitem{HW1}   T.Hibi, K.Watanabe, {\em   Study of three-dimensional algebras with straightening laws which are Gorenstein
domains. I}.  Hiroshima Math. J. 15 (1985), no. 1, 27--54. 

\bibitem{HW2}   T.Hibi, K.Watanabe,  { \em Study of three-dimensional algebras with straightening laws which are Gorenstein
domains. II}.  Hiroshima Math. J. 15 (1985), no. 2, 321--340.  

\bibitem{H2}   T.Hibi, {\em    Study of three-dimensional algebras with straightening laws which are Gorenstein domains. III.} 
Hiroshima Math. J. 18 (1988), no. 2, 299--308.

\bibitem{O} J. Oxley,  { \em Matroid Theory},  Oxford University Press, Oxford, 1992.  

\bibitem{R}  J.Roos,  {\em A description of the homological behaviour of families of quadratic forms in four variables}, in Syzygies
and Geometry, Boston 1995, A.Iarrobino, A.Martsinkovsky and J.Weyman eds., pp.86-95, Northeastern Univ. 1995. 

\bibitem{S}  B.Sturmfels,  {\em    Gr\"obner Bases and Convex Polytopes},  Amer. Math. Soc., Providence, RI,1995.

\bibitem{V} R.Villarreal, {\em Monomial algebras}. Monographs and Textbooks in Pure and Applied Mathematics, 238. Marcel Dekker,
Inc., New York, 2001. 

\bibitem{W1}    K.Watanabe, {\em Study of algebras with straightening laws of dimension $2$.} 
 Algebraic and topological theories (Kinosaki, 1984), 622--639, Kinokuniya, Tokyo, 1986.


\bibitem{W2}    K.Watanabe, {\em Study of four-dimensional Gorenstein ASL domains. I.
 Integral posets arising from triangulations of a $2$-sphere.}  Commutative algebra and combinatorics (Kyoto, 1985), 313--335, Adv.
Stud. Pure Math., 11, North-Holland, Amsterdam, 1987. 


\bibitem{We}   V.Weispfenning, {\em Comprehensive Gr\"obner bases}. J. Symbolic Comput. 14 (1992), no. 1, 1--29. 

 
\bibitem{Wel} D.Welsh, {\em Matroid Theory,} Academic Press,  London, 1976.


\bibitem{Wh} N. White,  { \em A unique exchange property for bases},  Linear Algebra Appl. 31 (1980), 81--91.
 
 \end{thebibliography}
\end{document}